\documentclass[11pt]{article}
\usepackage{ amsmath, amsthm, amssymb, amsfonts}
\usepackage[T1]{fontenc}
\usepackage[latin1]{inputenc}
\usepackage[extra]{tipa} 
\newcommand{\C}[1]{{\cal C}_{#1}}
\newcommand\inv{^{-1}} 
\newcommand\ov{\over}
\renewcommand\a{\alpha}
\renewcommand\b{\beta}
 \renewcommand\d{\delta}
 
\renewcommand{\s}{\sigma} 
\renewcommand\t{\tau}
\newcommand\ga{\gamma}

\newcommand\ve{\varepsilon}
\newcommand{\la}{\lambda}
 
\newcommand\cd{\cdots}
\newcommand\ld{\ldots}

\newcommand\iy{\infty}
\newcommand{\x}{\xi}

\newcommand{\noi}{\noindent}
 \def\({\left(} \def\){\right)}

\newcommand{\ra}{\rightarrow}
 
\newcommand{\bP}{\mathbb{P}}

\newcommand{\bZ}{\mathbb{Z}}

\newcommand{\cT}{\mathcal{T}}
\newcommand{\cC}{\mathcal{C}}

\newcommand{\cS}{\mathcal{S}}

\newcommand{\be}{\begin{equation}}
\newcommand{\ee}{\end{equation}}
\newcommand{\ba}{\begin{eqnarray*}}
\newcommand{\ea}{\end{eqnarray*}}
\newcommand{\bae}{\begin{eqnarray}}
\newcommand{\eae}{\end{eqnarray}}
\newcommand{\bc}{\begin{center}}
\newcommand{\ec}{\end{center}}

\newcommand{\fr}{\frac}

\begin{document}
\begin{titlepage}

\begin{center}{\Large\bf  Painlev\'e Functions in Statistical Physics}\end{center}
 
\begin{center}{\large Craig A.~Tracy} and {\large  Harold Widom}
\end{center}

\begin{abstract}
We review recent progress in limit laws for the one-dimensional asymmetric simple
exclusion process (ASEP) on the integer lattice.  The limit laws are expressed in terms of a certain
Painlev\'e II function.  Furthermore, we take this opportunity to give a brief survey of the appearance
of Painlev\'e functions in statistical physics.
\end{abstract}
\vspace{40ex}

\noindent Running Title: Painlev\'e in Statistical Physics
\par
\par
{\small
\begin{tabular}{ll}
\vspace{-1ex}
Craig A.~Tracy & Harold Widom\\
 \vspace{-1ex} Department of Mathematics& Department of Mathematics\\
 \vspace{-1ex}University of California&University of California\\
 \vspace{-1ex}Davis, CA 95616, USA&Santa Cruz, CA 95064, USA\\
 email: \texttt{tracy@math.ucdavis.edu}&email: \texttt{widom@ucsc.edu}
  \end{tabular} }

\end{titlepage}
  \begin{quote}
 ``It was a pleasant surprise to me
that such special functions actually appeared in
concrete problems of theoretical physics\ldots''
 \textit{Mikio Sato} \cite{andr}. 
 \end{quote}
\section{Introduction}
The appearance of Painlev\'e functions  in the 2D Ising model  is well-known \cite{MTW, WMTB}.  Equally well-known is that this provided one impetus for M.~Sato, T.~Miwa and M.~Jimbo \cite{SMJ} to develop their theory of \textit{holonomic quantum fields} which connects the theory of isomondromy preserving deformations of linear differential equations with the $n$-point correlation functions of the 2D Ising model.\footnote{A complete account of the SMJ theory can be found in the recent monograph by Palmer \cite{palmer}.}

The general consensus  in the field of ``exactly solvable models''  is that correlation functions
are expressible in terms of Painlev\'e functions  only in models
that are \textit{free fermion models}. More precisely, one expects  that for the appearance of functions of the Painlev\'e type, 
it is necessary for the underlying model or process to be a determinantal process in the sense of Soshnikov \cite{Sos2}.  In addition to the 2D Ising model, some notable  examples  where Painlev\'e functions arise in correlation functions include the
one-dimensional impenetrable Bose gas \cite{FFGW, JMMS, Koj, KBI}, the  
Ising chain in a tranverse field \cite{PCQN}, 
the distribution functions of random matrix theory \cite{AvM, Baik, BD, FW, JMMS, TW1, TW2, TW3}, Hammersley's growth process \cite{BDJ, BR1},
corner and polynuclear growth models \cite{BR2, GTW, Joh1, PS1, PS2} and
the totally asymmetric simple exclusion process (TASEP) \cite{BAC, Joh1, PS3}.  Universality theorems in random matrix theory have  extended the appearance of Painlev\'e functions to a wide class of matrix ensembles \cite{BI, DG1, DG2, DKMVZ,  Sos1}.\footnote{It is also worth noting 
that due to the close connection of random matrix theory to multivariate statistical analysis, these same distribution functions involving Painlev\'e functions are now routinely used in data analysis \cite{John1, John2, PPR}.}  
In recent developments \cite{ACQ, SS1, SS2, SS3}  Painlev\'e II appears in the long time asymptotics
of explicit  formulas
for the exact height distribution for the KPZ equation \cite{KPZ} with narrow wedge initial condition.

  As just noted, one does not expect Painlev\'e functions to arise in correlation functions in models that are exactly solvable in the sense of Baxter \cite{Bax} but are not free fermion models, e.g.\ 6-vertex model, XXZ quantum spin chain, Baxter's 8-vertex model.
Having said that, the \textit{universality conjecture} arising in the theory of
phase transitions  suggests, for instance, that the scaling limit of a large class of ferromagnetic 2D Ising models is the same as that of the Onsager 2D Ising model; and hence, Painlev\'e functions are conjectured to appear 
(in the massive scaling limit)  in models outside of the class of exactly solvable models.  This last statement is substantiated by the developments
in \cite{ACQ, SS1, SS2, SS3}. 

In this paper we review recent progress \cite{TW4, TW5, TW6, TW7, TW8} on the current fluctuations in the asymmetric simple exclusion process (ASEP) on the integer lattice $\bZ$ \cite{Li1, Li2}.  ASEP is in the class of Bethe Ansatz solvable models \cite{deGE, GS}
 but only for certain values of the parameters is ASEP a
determinantal process \cite{Joh1, PS3, Sch1}. 
That ASEP is Bethe Ansatz solvable comes as no surprise \textit{once} one realizes that the generator of ASEP
is a similarity (not unitary!) transformation of the XXZ-quantum spin Hamiltonian \cite{ADHR, Sch2, YY}.
Our main results relate the limiting current fluctuations in ASEP for certain initial conditions to the TW distributions $F_1$ and $F_2$ of random matrix theory \cite{TW3, TW4}.  Both $F_1$ and $F_2$  are expressible in terms of the same Hastings-McLeod solution of Painlev\'e II  \cite{FIKN, HMc}, see \S4.2. 

\section{Master Equation and Bethe Ansatz Solution}
Since its introduction in 1970 by F.~Spitzer \cite{Spi}, the asymmetric simple exclusion
process (ASEP) has attracted considerable attention both in the mathematics and physics literature due to the fact
it is one of the simplest lattice models describing transport far from equilibrium.
Recall \cite{Li1, Li2} that the ASEP on the integer lattice $\bZ$ is a continuous time
 Markov process $\eta_t$
where $\eta_t(x)=1$ if $x\in \bZ$ is occupied at time $t$, and $\eta_t(x)=0$ if $x$ is vacant at time $t$.
Particles move on $\bZ$ according to two rules: (1) A particle at $x$ waits an exponential time with
parameter one, and then chooses $y$ with probability $p(x,y)$;
 (2) If $y$ is vacant at that time it moves
to $y$, while if $y$ is occupied it remains at $x$.  The adjective ``simple'' refers to the fact that the allowed
jumps are only one step to the right, $p(x,x+1)=p$, or one step to the left, $p(x,x-1)=q=1-p$.  The
totally asymmetric simple exclusion process (TASEP) allows jumps only to the right ($p=1$) or only to the left ($p=0$).\footnote{It is TASEP that is a determinantal process.} In the mapping from the XXZ quantum spin chain, the anisotropy parameter $\Delta$ of the spin chain is related to the hopping probabilities $p$ and $q$ by
\[ \Delta=\fr{1}{2\sqrt{pq}}\ge 1,\]
the ferromagnetic regime of the XXZ spin chain.

We begin with a system of $N$ particles and later take the limit $N\ra\iy$.  A configuration is specified by giving
the location of the $N$ particles. We denote by $Y=\{y_1,\ldots,y_N\}$ with $y_1<\cd<y_N$ the initial configuration of the process and write $X=\{x_1,\ldots,x_N\}\in\bZ^N$. When $x_1<\cdots<x_N$ then $X$ represents a possible configuration of the system at a later time $t$. We denote by $P_Y(X;t)$ the probability that the system is in configuration $X$ at time $t$, given that it was initially in configuration $Y$.

Given $X=\{x_1,\ld,x_N\}\in\bZ^N$ we set
\[X_i^+=\{x_1,\ldots,x_{i-1},x_i+1,x_{i+1},\ldots,x_N\},\quad
X_i^-=\{x_1,\ldots,x_{i-1},x_i-1,x_{i+1},\ldots,x_N\}.\]
The master equation for a function $u$ on $\bZ^N\times\mathbb{R}^+$ is
\be{d\ov dt}u(X;t)=\sum_{i=1}^N\Big(p\,u(X_i^-;t)+q\,u(X_i^+;t)-u(X;t)\Big),\label{master}\ee
and the boundary conditions are, for $i=1,\ld,N-1$,
\[u(x_1,\ldots,x_i,x_i+1,\ldots,x_N;t)\]
\be = p\, u(x_1,\ld, x_i,x_i,\ld,x_N;t)+q\, u(x_1,\ld, x_i+1,x_i+1,\ldots,x_N;t).\label{boundary}\ee
The initial condition is
\be u(X;0)=\d_Y(X)\ \ {\rm when}\ x_1<\cdots<x_N.\label{initial} \ee
The basic fact is that if $u(X;t)$ satisfies the master equation, the boundary conditions, and the initial condition, then $P_Y(X;t)=u(X;t)$
when $x_1<\cd<x_N$.  This is, of course, one of Bethe's basic ideas (see, e.g., \cite{Bat}):  incorporate the interaction (in this case
the exclusion property)  into the boundary conditions (\ref{boundary}) of a free particle system (\ref{master}).

Recall that an inversion in a permutation $\s$ is an ordered pair $\{\s(i),\s(j)\}$ in which $i<j$ and $\s(i)>\s(j)$. We define \cite{YY}
\be S_{\a\b}=-{p+q\x_\a\x_\b-\x_\a\ov p+q\x_\a\x_\b-\x_\b}\label{Smatrix}\ee
and then
\[A_\s=\prod\{S_{\a\b}:\{\a,\b\}\ {\rm is\ an\ inversion\ in}\ \s\}.\]
We also set
\[\ve(\x)=p\,\x\inv+q\,\x-1.\]  In the next theorem we shall assume $p\ne 0$, so the $A_\s$ are analytic at zero in all the variables. Here and later all differentials $d\x$ incorporate the factor $(2\pi i)\inv$. 

\noi{\bf Theorem 2.1}. We have
\be P_Y(X;t)=\sum_{\s\in\cS_N}\int_{\C{r}}\cdots\int_{\C{r}} A_\s\,\prod_i\x_{\s(i)}^{x_i-y_{\s(i)}-1}\,e^{\,\sum_i\ve(\x_i)\,t}\,d\x_1\cd d\x_N,\label{theorem1}\ee
where $\C{r}$ is a circle centered at zero with radius $r$ so small that all the poles of the integrand lie outside $\C{r}$.
 
The proof that $P_Y(X;t)$ satisfies (\ref{master}) is immediate and the fact it satisfies the boundary conditions (\ref{boundary}) is exactly the same argument as in the XXZ problem \cite{YY}.  The difficulty lies in showing (\ref{theorem1}) satisfies the initial condition (\ref{initial}).  Observe that the term in (\ref{theorem1}) corresponding to the identity permutation does satisfy the initial condition.  Thus the proof will be complete once one demonstrates that the remaining $n!-1$ other terms sum to zero at $t=0$.  This is indeed the case (some are individually zero and others cancel in pairs) and the result depends crucially upon the choice of the contours $\C{r}$ \cite{TW4}.  For the special case of TASEP,  $p=1$,  it follows from (\ref{Smatrix}) and (\ref{theorem1}) that
the right-hand side of (\ref{theorem1}) can be expressed as a $N\times N$ determinant as first obtained in \cite{Sch1}.

We note that unlike the usual applications of Bethe Ansatz, it is not the spectral theory of the operator that is of interest but rather the
transition probability $P_Y(X;t)$.  Thus there are no Bethe equations in our approach; and hence, no issues concerning the completeness of the Bethe eigenfunctions.  Indeed, there is not even an Ansatz in this approach!  We remark that this result  extends with only
minor modifications  to the solution 
$\Psi(x_1,\ldots, x_N;t)$ of the time-dependent Schr\"odinger equation with XXZ Hamiltonian where the $x_i$'s denote the location of the $N$  ``up spins'' in a sea of ``down spins'' on $\bZ$.  

\section{Marginal Distributions and the Large $N$ Limit}

We henceforth assume $q>p$ so there is a net
drift of particles to the left.    
Here we consider two different initial conditions.  The first, called \textit{step initial condition}, starts with particles located at $\bZ^+=\{1,2,\ldots\}$.   The second initial
condition is the 
\textit{step Bernoulli initial condition}: each site in $\bZ^+$, independently of the others, is initially occupied with probability $\rho$,
$0<\rho\le 1$; all other sites are initially unoccupied.  In each of these cases it makes sense to speak of the position of the $m$th particle from the left at time $t$, $x_m(t)$, and its distribution function $\bP(x_m(t)\le x)$.  It is elementary to relate $\bP(x_m(t)\le x)$ to the distribution of the
total current $\cT$ at position
$x$ at time $t$,
\[ \cT(x,t):=\textrm{number of particles}\> \le x\> \textrm{at time}\>\> t;\]
namely,
\[ \bP(\cT(x,t)\le m)=1-\bP(x_{m+1}(t)\le x).\]
For this reason we first concentrate on $\bP_Y(x_m(t)\le x)$ and only at the end translate the results into statements concerning $\cT$. 
(The subscript $Y$ denotes the initial configuration.)

Now for finite $Y$
\[ \bP_Y(x_m(t)= x)= \sum_{x_1<\cdots <x_{m-1}<x<x_{m+1}<\cdots <x_N} P_Y(x_1,\ldots, x_{m-1},x,x_{m+1},\ldots, x_N;t).\]
Since the contours $\C{r}$  in (\ref{theorem1}) have $r\ll 1$, the sums over $x_{m+1},\ldots, x_N$ can be interchanged with the integrations 
in variables $\x_{\s(j)}^{x_{j}}$, $m+1\le j\le N$, and the geometric series summed. To perform the sums over $x_1, \ldots, x_{m-1}$ the contours
in the $\x_{\s(j)}^{x_j}$ variables, $1\le j\le m-1$,  must be deformed out beyond the unit circle and then the sums can be
interchanged with the integrations.  This deformation
beyond the unit circle can be done in such a way as not to encounter any poles of the integrand.
  However, upon deforming these contours back to $\C{r}$ (after the geometric series are summed) one does encounter poles; and one finds
some remarkable cancellations: only the residues from the poles at $\x_i=1$  are nonzero. The result is a sum over all subsets of $S$ of
$\{1,\ldots, N\}$ with $|S^c|<m$ whose summands involve $|S|$-dimensional integrals with contours $\C{r}$.\footnote{This is Theorem 5.1 in \cite{TW4}.}
However, this resulting expression for $\bP_Y(x_m(t)=x)$ is not so useful for taking the $N\ra\iy$ limit.

The next step is to expand the contours to $\C{R}$, $R\gg 1$.  It is then possible to take the $N\ra \iy$ limit in the resulting expressions.  
The details \cite{TW4} are involved and they depend crucially upon some algebraic identities which we now state.

\subsection{Three identities}
Let
\[ f(i,j):=p+q\x_i\x_j-\x_i.\]
Identity \#1:
\be
\sum_{\s\in\cS_N} \textrm{sgn}(\s) \,
\fr{\prod_{i<j} f(\s(i),\s(j))}{(\x_{\s(1)}-1)(\x_{\s(1)}\x_{\s(2)}-1)\cdots (\x_{\s(1)}\x_{\s(2)}\cdots \x_{\s(N)}-1)}=
q^{N(N-1)/2}\,\fr{\prod_{i<j}(\x_j-\x_i)}{\prod_j(\x_j-1)}.
\label{iden1}\ee
Identity \#2: For $N\ge m+1$,
\be
\sum_{|S|=m} \prod_{{i\in S\atop j\in S^c}}\fr{f(i,j)}{\x_j-\x_i}\, \left(1-\prod_{j\in S^c} \x_j\right) = q^m \left[{N-1 \atop m}\right] (1-\prod_{j=1}^N\x_j)
\label{iden2}\ee
In (\ref{iden2}) the sum runs over all subsets $S$ of $\{1,\ldots, N\}$ with cardinality $m$, and $S^c$ denotes the complement of $S$ in
$\{1,\ldots, N\}$. Here $\left[{N\atop m}\right]$ is a slightly modified  $\tau$-binomial coefficient, $\tau:=p/q$,
\ba
[N]  &:=& \fr{p^N-q^N}{p-q}, \>\> [0]:=1,\\
\left[N\right]! &:=& [N] \,[N-1] \,\cdots  [1],\\
\left[{N\atop m}\right]&:=& \fr{[N]!}{[m]! [N-m]!}=q^{m(N-m)}\left[{N\atop m}\right]_\tau\ea
where $\left[{N\atop m}\right]_\tau$ is the usual $\tau$-binomial coefficient.  We define $\left[{N\atop m}\right]_\tau=0$ for $m<0$.
In proving (\ref{iden2}) we first proved the simpler identity

\noi Identity \#3:
\[ \sum_{|S|=m} \prod_{{i\in S\atop j\in S^c}}\fr{f(i,j)}{\x_j-\x_i} = \left[{N \atop m}\right]. \]
We believe that these identities suggest a deeper mathematical structure that is yet to be discovered.

\subsection{Final expression for $\bP(x_m(t)\le x)$ for step and step Bernoulli initial conditions}
We denote by $\bP_\rho$ the probability measure for ASEP with step Bernoulli initial conditions.  For $\rho=1$ the measure is
ASEP with step initial condition.  Let
\[ c_{m,k}:= (-1)^{m}\, q^{k(k-1)}  \tau^{m(m-1)/2} \tau^{-km} \left[{k-1\atop m-1}\right]_\tau .\]
Observe that $c_{m,k}=0$  when $m>k$.

\noi\textbf{Theorem 3.1} \cite{TW4, TW8}.
Assume $q>p$, then
\bae \bP_\rho(x_m(t)\le x)&=&\sum_{k\ge 1} \fr{q^{k(k-1)/2}\t^{k(k+1)/2}}{k!}\,  c_{m,k}  \int_{\C{R}}\cdots \int_{\C{R}} \prod_{1\le i\neq j\le m} \fr{\x_j-\x_i}{f(i,j)}
\times \nonumber \\ 
&& \prod_i\fr{\rho}{\x_i-1+\rho(1-\t)}\,\prod_{i =1}^m\fr{\x_i^{x} e^{t\ve(\x_i)}}{1-\x_i}\, d\x_i 
\label{marginal}\eae
The contour $\C{R}$, a circle of radius $R\gg 1$ centered at the origin,  is chosen so that all (finite) poles of the integrand lie inside the contour.

We remark that for TASEP, $p=0$, the above sum reduces to one term; and this term can be shown to be equal to a $m\times m$ determinant.

The final simplification results if we use the identity \cite{TW5}
\[ \det\left(\fr{1}{f(i,j)}\right)_{1\le i,j\le k}=(-1)^k (pq)^{k(k-1)/2} \prod_{i\neq j}\fr{(\x_j-\x_i)}{f(i,j)}\, \prod_i\fr{1}{(1-\x_i)(q\x_i-p)}\]
in (\ref{marginal}) and recognize the summand, a $k$-dimensional integral, as the coefficient of $\la^k$ in the Fredholm expansion
of $\det(I-\la K_\rho)$ where
$K_\rho$ acts on functions on $\C{R}$ by
\[ f(\x)\longrightarrow \int_{\C{R}} K_\rho(\x,\x') f(\x')\,d\x' \]
where
\be K_\rho(\x,\x')= q \,\fr{\x^x e^{t\ve(\x)}}{p+q\x\x'-\x}\, \fr{\rho(\x-\t)}{\x-1+\rho(1-\t)}, \>\>\t=\fr{p}{q}.\label{kernel}\ee
Note that when $\rho=1$, the case of step initial condition, the last factor in $K_\rho(\x,\x')$ equals one.

Since the coefficient of $\la^k$ in the expansion of $\det(I-\la K_\rho)$ is equal to
\[\fr{(-1)^k}{k!}\int \det(I-\la K_\rho)\, \fr{d\la}{\la^{k+1}},\]
this fact together with the $\t$-binomial theorem  gives the final result for $\bP_\rho(x_m(t)\le x)$.

\noi\textbf{Theorem 3.2} \cite{TW4, TW8}.  Let $\bP_\rho$ denote the probability measure for ASEP with step Bernoulli initial condition with density
$\rho$ and   $x_m(t)$ denote the position of the $m$th particle from the left at time $t$, then
\be \bP_\rho(x_m(t)\le x)= \int_{\cC} \fr{\det(I-\la K_\rho)}{\prod_{j=0}^{m-1} (1-\la\t^j)} \, \fr{d\la}{\la} \label{marginal2}\ee
where the contour $\cC$ is a circle centered at the origin enclosing all the singularities at $\la=\t^{-j}$, $0\le j\le m-1$ and $K_\rho$ is the integral
operator whose kernel is given by (\ref{kernel}).

\section{Limit Theorems}
\subsection{KPZ Scaling}
The scaling limit that is of most interest is the \textit{KPZ scaling limit} \cite{KPZ, Spo}.  In the terminology here this scaling limit is
\[ m\ra\iy, \> t\ra\iy\>\> \textrm{with}\>\> \s=\fr{m}{t}\le 1\>\>\textrm{fixed}. \]
As we shall see, the limiting distribution will depend upon the relative sizes of $\s$ and $\rho^2$.  For the moment we concentrate on the cases
$0<\s < \rho^2$ and $\s=\rho^2$ with $0<\rho\le 1$.
As in any central limit theorem, to obtain a nontrivial limit the $x$ in $\bP_\rho(x_m(t)\le x)$ must also be scaled (this too is part
of KPZ scaling).   In anticipation of
the theorem we set
\[ x:= c_1 t + c_2 t^{1/3} s \]
where the $\tfrac{1}{3}$ is the famous KPZ universality exponent \cite{KPZ, MMMT} and 
\[ c_1:=-1+2\sqrt{\s}, \>\> c_2:=\s^{-1/6}(1-\sqrt{\s})^{2/3}.\]
The two distribution functions that arise in the KPZ scaling limit are defined in the next section.

\subsection{Distributions $F_1$ and $F_2$}
The distributions $F_1$ and $F_2$ can be defined by either their Fredholm determinant representations or their representations in terms of a Painlev\'e II
function.  Here we take the latter route.
Let $q$ denote the solution to the  Painlev\'e II equation
\[q''=x \,q + 2 q^3\]
satisfying
\[ q(x) \sim \textrm{Ai}(x), \> x\ra\iy,\]
where $\textrm{Ai}(x)$ is the Airy function.  That such a solution exists and is unique was proved by Hastings and McCleod \cite{HMc}.\footnote{A modern account of Painlev\'e transcendents can be found in the monograph by Fokas, et al.\ \cite{FIKN}.}
Then we have
\bae F_2(s)&=& \exp\left(-\int_s^\iy (x-s) q^2(x)\, dx\right),\label{F2}\\
F_1(s)&=& \exp\left(-\fr{1}{2}\int_s^\iy q(x)\, dx\right) \left(F_2(s)\right)^{1/2}. \label{F1}
\eae
The asymptotics of these distributions as $x\ra\iy$ is straightforward given the large $x$ asymptotics of the Airy function; however, the complete
asymptotic expansion as $x\ra -\iy$ has only recently been completed \cite{BBD}.  For high-accuracy numerical evaluation of $F_1$ and $F_2$,
it turns out that it is better to start with their Fredholm determinant representations \cite{Born}.

\subsection{Limit Laws}
The asymptotic analysis \cite{TW6, TW8}
 of the Fredholm determinant in the formula for  $\bP_\rho(x_m(t)\le x)$ in (\ref{marginal2})
required the development of new methods since the operator $K_\rho$ is not of the usual ``integrable integral operator'' form normally appearing
in random matrix theory \cite{Blo, IIKS, TW2}.  The main point is that the kernel $K_\rho$ has the same Fredholm determinant as sum of two kernels;
one has large norm but fixed spectrum and its resolvent can be computed exactly, and the other is better behaved \cite{TW6}.

We now state the results of this asymptotic analysis.

\noi\textbf{Theorem 4.1} \cite{TW6, TW8}.  When $0\le p< q$, $\ga:=q-p$, 
\bae
\lim_{t\ra\iy}\bP_\rho\left({x_m(t/\ga)-c_1 t\ov c_2 t^{1/3}}\le s\right)= F_2(s)\>\>\textrm{when}\>\> 0<\s<\rho^2,\label{limit1}\\
\lim_{t\ra\iy}\bP_\rho\left({x_m(t/\ga)-c_1 t\ov c_2 t^{1/3}}\le s\right)= F_1(s)^2\>\>\textrm{when}\>\> \s=\rho^2, \>\rho<1.\label{limit2}
\eae

\noi This theorem implies a limit law for the current fluctuations.  Define
\[ v=x/t, \>\> a_1=(1+v)^2/4,\>\> a_2=2^{-4/3}(1-v^2)^{2/3}.\]

\noi\textbf{Theorem 4.2}. When $0\le p< q$, $\ga:=q-p$, 
\bae
\lim_{t\ra\iy}\bP_\rho\left({\cT(vt,t/\ga)-a_1 t\ov a_2 t^{1/3}}\le s\right)= 1-F_2(-s)\>\>\textrm{when}\>\> -1<v<2\rho-1,\label{limit3}\\
\lim_{t\ra\iy}\bP_\rho\left({\cT(vt,t/\ga)-a_1 t\ov a_2 t^{1/3}}\le s\right)= 1-F_1(-s)^2\>\>\textrm{when}\>\> v=2\rho-1, \>\rho<1.\label{limit4}
\eae

\noi For step initial condition with $0<\s<1$ the limit laws are (\ref{limit1}) and (\ref{limit3}) \cite{TW6, TW7}.  When $\s>\rho^2$ (or $v>2\rho-1$)
the fluctuations are of order $t^{1/2}$ and the limiting distribution is Gaussian, see \cite{TW8} for details.

For TASEP, $p=0$,  with step initial condition the limit law (\ref{limit3})  was first proved by Johansson \cite{Joh1}.  
For TASEP with step Bernoulli initial condition the limit laws (\ref{limit3}) and (\ref{limit4})
were conjectured by Pr\"ahofer and Spohn \cite{PS3} and proved recently by Ben Arous and Corwin \cite{BAC}.  The fact that these limit laws remain essentially identical (the only change is the factor $\ga$ in the time slot) is a very strong statement of KPZ Universality.  From the integrable systems perspective, these results are, to the best of the authors' knowledge, the first limit laws of Bethe ansatz solvable models (outside the class of determinantal models) where the correlation functions are expressible in terms of Painlev\'e functions.

\begin{table}
\begin{center}
\caption{The mean ($\mu_\beta$),  variance ($\sigma^2_\beta$),
skewness ($S_\beta$) and  kurtosis ($K_\beta$) of $F_\beta$, $\beta=1,2$.
The  numbers are courtesy of F.~Bornemann and M.~Pr\"ahofer.} 
\vspace{4ex}
\begin{tabular}{|l|llll|}\hline\vspace{2pt}
$\beta$ & \qquad$\mu_\beta$ & \qquad$\sigma^2_\beta$ & \qquad$S_\beta$ & \qquad$K_\beta$ \\  \hline
1 & -1.206\,533\,574\,582 & 1.607\,781\,034\,581 & 0.293\,464\,524\,08 & 0.165\,242\,9384 \\
2 & -1.771\,086\,807\,411 & 0.813\,194\,792\,8329 & 0.224\,084\,203\,610 & 0.093\,448\,0876 \\
 \hline
\end{tabular}
\end{center}
\end{table}

\section{Conclusions}
Today Painlev\'e functions  occur in many areas of theoretical statistical physics.  In the case
of KPZ fluctuations there are now experiments \cite{MMMT, Tak} on stochastically growing interfaces where quantities such as  the skewness and the kurtosis of $F_\beta$  (see Table 1), as well as
the distribution functions themselves,  are compared with experiment.   In \cite{Tak}  K.~Takeuchi and M.~Sano conclude that
their measurements    ``\ldots have shown without
fitting that the fluctuations of the cluster local radius
asymptotically obey the Tracy-Widom distribution of the
GUE random matrices.''

\noindent{\textsc{Acknowledgements:}}  This work was supported
by the National Science Foundation under grants DMS--0906387 (first author)
and DMS--0854934 (second author).

\end{document}